\documentclass[11pt]{article}
\usepackage{amsmath,amsfonts,enumerate}

\def\R{{\mathbb R}}
\def\Z{{\mathbb Z}}
\def\Q{{\mathbb Q}}
\def\C{{\mathbb C}}
\def\H{{\mathbb H}}
\def\rank{{\mathrm{rank}\,}}
\newtheorem{theorem}{Theorem}
\newtheorem{lemma}{Lemma}
\newtheorem{corollary}{Corollary}

\begin{document}

\title{The type numbers of closed geodesics
\thanks{A translation of an appendix to the Russian edition of
``Calculus of variations in the large'' by M. Morse.
The work was supported by RFBR (grant 09-01-00598) and
Max Planck Institute for Mathematics in Bonn.}}
\author{I.A. Taimanov
\thanks{Institute of Mathematics, 630090 Novosibirsk, Russia;
e-mail: taimanov@math.nsc.ru.}
}
\date{}
\maketitle


The seminal book \cite{Morse} by M. Morse  had put (together with
the treatise by Lus\-ter\-nik and Schnirelmann \cite{LS})
the foundations for applying
topological methods to variational problems.

Writing this appendix to its Russian translation we had no intention
to survey wide applications of the Morse theory but wanted to
complement the text by clarifications of some statements related to
the type numbers of closed geodesics.

The foundations of what is called now the Morse theory together with
their first applications are exposed in Chapters 1--7.

The last chapters 8 and 9 are concerned with applications of the
Morse theory to the theory of closed geodesics.

In particular, the type numbers of closed geodesics in
the space $\Omega$ of closed non-oriented curves without marked points
are studied in Chapter 8. In
\cite{Schwarz} it was noted that the statement about these type
numbers is not correct for iterated geodesics
and the necessary correction is given.

In Chapter 9 this statement on the type numbers is used for computing
the relative homology groups $H_\ast(\Omega(S^n),S^n; \Z_2)$ for spheres
(here $S^n \subset \Omega(S^n)$ is the set formed by one-point curves).
Due to the mistake concerning the type numbers
these computations are also wrong and these homologies are
not computed until recently. However,
the main problem which is studied in
Chapter 9 is the existence of closed geodesics on manifolds sufficiently closed
to ellipsoids. One of the main result proved by Morse
is as follows:

\begin{itemize}
\item
{\sl on an ellipsoid in $\R^{n+1}$ which is close to the round sphere but has
pairwise different axes every prime closed geodesic is either one of
$\frac{n(n+1)}{2}$ principal ellipses, either is very long.}
\end{itemize}

The short modern proof of this theorem is given in \cite{KlingenbergII}.
In the 1980s that was generalized by Ballmann \cite{Ballmann83}
(for $n=2$) and then by Bangert
\cite{Bangert86} (for all $n$) who proved that

\begin{itemize}
\item
{\sl for every $\varepsilon > 0$
there exists $\delta=\delta(\varepsilon,n)$ such
that every prime closed geodesic on a Riemannian sphere $S^n$ whose sectional
curvature $K$ satisfies $1-\delta < K < 1+\delta$ has length in
$(2\pi-\varepsilon,2\pi+\varepsilon)$ or larger than
$\varepsilon^{-1}$.}
\end{itemize}

Other results by Morse did also concern closed geodesics on spheres
and were also formulated for metrics sufficiently closed to the
round metrics. The modern approach to this problem was proposed by
Ballmann, Thorbergsson, and Ziller who used the comparison theory
and established the existence of simple closed geodesics for spheres
with $\frac{1}{4}$-pinched metrics \cite{BTZ83}.

As we already mentioned the groups
$H_\ast(\Omega(S^n),S^n; \Z_2)$
are still not known however some nontrivial cycles were
detected and were used for proving the existence of closed geodesics.
This story did also involve some mistakes and
``non-existing'' cycles. The most recent description of
these homologies is given in \cite{Anosov81,Anosov82}.

The structure of our text is as follows.

In \S\S 1--3 we give a short introduction into the type numbers of closed
geodesics and their applications and expose the necessary
corrections to Morse'e results.

In \S 4 we demonstrate how these corrections were applied. Therewith
we present the Bott formula for indices of iterated geodesics which
is necessary for the demonstration and the Gromoll--Meyer theorem
which is derived from the Bott formula.

Since all preliminary material
is already done we could not resist a temptation to give in \S 5 a short
survey of main methods and results concerning the application of the
Morse theory \cite{Morse} to proving the existence of
closed geodesics.

In \S 6 we briefly remark on two active directions
of developing the Morse theory and hope that this section shall invite
researchers into the subject.

We thank D.V. Anosov, W. Ballmann, and H.-B. Rademacher
for helpful discussions.

\section{The type numbers of closed geodesics}

In terms of relative homology groups the definition of type
numbers \cite[p.222]{Morse} is as
follows \cite{ST}.

Let
$$
F: X \to \R
$$
be a smooth function. For arbitrary $r$ we put
$$
X_r = \{F < r\} \subset X, \ \ \ \bar{X}_r = \{F \leq r\}.
$$
Let $c$ be a critical level of $f$ and $\omega \subset F^{-1}(c)$ be
a set of critical points which probably does not contain all
critical points with $F=c$. Let $U$ be a neighborhood of $\omega$
which does not contain other critical points than from $\omega$ and
which is contractible to $\omega$ and let $\varepsilon > 0$ be a
positive constant such that the set
$$
X_c \setminus \bar{X}_{c-\varepsilon} = \{c-\varepsilon \leq F < c\}
$$
contains no critical points of $F$. Then the $k$-th type number
$m_k$ of $\omega$ is the $k$-th Betti number of the pair
$(X_{c-\varepsilon} \cup U, X_{c-\varepsilon})$:
$$
m_k (\omega) = \rank H_k(X_{c-\varepsilon} \cup U, X_{c-\varepsilon}).
$$
This definition precedes Theorem 10.1 which is as follows:

\begin{itemize}
\item
{\sl let $M_k$ be the sum of the $k$-th type numbers of all critical
sets in $X_a$ and $B_k$ be the $k$-th Betti number of $X_a$, then
$$
M_0 \geq B_0,
$$
\begin{equation}
\label{morseineq} M_ 0 - M_1 \leq B_0 - B_1,
\end{equation}
$$
\dots,
$$
$$
M_0 - M_1 + \dots + (-1)^r M_r =
B_0 - B_1 + \dots + (-1)^r B_r
$$
for any sufficiently large integer $r$.} \footnote{Here we assume
that $X$ is compact or that $X_a$ contains finitely many connected
critical sets $\omega$ for which all type numbers are finite and all
of them except finitely many vanish.}
\end{itemize}

The inequalities (\ref{morseineq}) are well-known as
{\it the Morse inequalities} and, in particular, they imply that
$$
M_k \geq B_k,
$$
thus enabling us to prove the existence of critical points which are
``{\it topologically necessary}''.
\footnote{The Betti numbers are called the connectivities in \cite{Morse}.}
The definition of the type numbers as well as of the Betti numbers depends
on the choice of the coefficient group. Morse considered
homologies with coefficients in $\Z$ and $\Z_2$.

Closed geodesics on a Riemannian manifold $M$ are critical points of
the length functional on the admissible spaces of  closed curves:
\footnote{\label{space} Since the analytical aspects is not our main
concern we avoid here rigorous definitions due to analytical
complications. In fact, the length functional is not differentiable
even on a smooth parameterized curve $r(t)$ such that $\dot{r} = 0$
at some point. Hence, usually one have not only to restrict
ourselves to curves for which the length is well-defined but to pass
from spaces of such curves to their deformation retracts in which
every closed geodesic line is represented and on which the length is
defined as a smooth functional. This approach was applied by Morse.}
$$
r: [0,1] \to M, \ \ \ r(0)=r(1).
$$
On this space acts the orientation-reversing map
$$
\sigma(r) = s \ \ \ \mbox{if $r(t) =s(1-t)$ for all $t$}.
$$
Since the length is independent of a parametrization we have to
identify curves $r$ and $s$ if they are related by a
change of parametrization:
$$
r \sim s \ \ \ \mbox{iff} \ \ \
r(t) = s(f(t)) \ \ \mbox{with $f:[0,1] \to [0,1]$ and $f^\prime \geq 0$}.
$$
Therewith we obtain {\it the space $P^+(M)$ of oriented closed
curves with a marked point} which is $r(0)$. Any curve from $P^+(M)$
is represented uniquely by a map $r: [0,1] \to M$ with the length
parameter $t$, i.e. $|\dot{r}|=\mathrm{const}$.

The quotient space
$$
P(M) = P^+(M)/\sigma
$$
is {\it the space of non-oriented closed curves with a marked
point.}

On the space $P^+(M)$
there is a group action
\begin{equation}
\label{theta}
g_\theta: P^+(M) \to P^+(M), \ \ \theta \in S^1 =
\R/\Z,
\end{equation}
of $S^1 = SO(2)$ which changes the
marked points.
Let $\gamma \in P^+(M)$ is represented by
a map
$r: [0,1] \to M$ with the length parameter $t$.
Then
$g_\theta(\gamma)$ is represented by
$$
r_\theta(t) = r((t+\theta)\mathrm{mod}\,1).
$$
Since this action transforms closed closed into closed geodesics, to
every nontrivial (i.e., not one-point curve) closed geodesic
there corresponds a circle of critical points of the length functional.

By the involution $\sigma$ the action of $SO(2)$ is extended to
the action of $O(2)$.

{\it The periodic problem of the calculus of variations} consists in
estimating from below the number of closed extremals of functionals.
In Riemannian geometry these extremals are geodesics.

To apply topological methods to finding
closed geodesics we have to evaluate the number of critical points of
the length functional on the admissible space of closed curves.
Here we come to two problems:

1) in the space of non-parameterized curves with marked points to
every closed geodesic line there corresponds an $S^1$-family of
critical points;

2) in the quotient spaces
$$
\Omega^+(M) = P^+(M)/SO(2),
$$
{\it the space of
oriented closed curves} (without marked points),
and
$$
\Omega(M) = P^+(M)/O(2),
$$
to every nontrivial
closed geodesic $\gamma$ there corresponds infinitely many critical
points given by iterations of $\gamma$. So we say that
\footnote{
It is clear that these definitions are applicable to closed
extremals of all functionals defined by some Lagrangian functions:
$$
F(\gamma) = \int_\gamma L(x,\dot{x})\, dt,
$$
since for them iterations of extremals are also extremals.}

\begin{itemize}
\item
a closed geodesic is {\it prime} if it is not
a one-point geodesic or an iteration of another closed geodesic;

\item
a prime geodesic is called {\it simple} if it has no selfintersections;

\item
two closed geodesics are {\it geometrically distinct}
if the corresponding geodesic lines on the manifold are different;

\item
two closed geodesics are called {\it algebraically distinct} if
either they are geometrically distinct, either they are different
iterations of the same prime closed geodesic.
\end{itemize}

Therefore different critical points of the length functional on
$$
F: \Omega^+ = P^+/S^1 \to \R
$$
correspond to algebraically distinct closed geodesics.
The way to prove the existence of closed extremals is as follows:

1) to estimate the Betti numbers of $\Omega^+(M)$ in terms of the topology
of $M$;

2) derive the existence of ``topologically necessary''
critical points of the length functional on
$\Omega^+(M)$ from the Morse inequalities (\ref{morseineq}).

\medskip

{\sc Remark.}
Another approach to look for closed geodesics, which is
the most adopted now, is to consider the energy functional
on the spaces of absolutely continuous closed curves with finite energy:
$$
\Lambda M = \{ r: S^1 \to M, \ \mbox{$r$ is absolutely continuous},
\int |\dot{r}|^2 \, dt < \infty\}
$$
(see \cite{Klingenberg}). This is the so-called {\it free loop space}
and it is a Hilbert manifold.
Such an approach has many nice features:

1) extremals of the energy functional on $\Lambda M$ are exactly
length-para\-me\-te\-ri\-zed geodesics;

2) any smooth
map $f:M \to N$ induces a smooth
map $\Lambda f: \Lambda M \to \Lambda N$ (if $f$ is only Lipschitz
then $\Lambda f$ may be discontinuous):

3) the action of $O(2)$ is defined on $\Lambda M$
in the standard way and there is an equivariant continuous map
$\Lambda M \to \Lambda M$ which is homotopical to the identity and maps
$\Lambda M$ onto the subspace of length-parameterized curves. This implies that
$\Lambda M, \Lambda M/SO(2)$ and their $\sigma$-quotients
have the same homotopical properties as $P^+(M), \Omega^+(M)$ and their
$\sigma$-quotients.

We refer to \cite{Anosov80} for proofs of the statements 2) and 3).

Since the action of $SO(2)$ is not free (it branches exactly at iterated
curves), $\Lambda M/SO(2)$ is ``a Hilbert
orbifold'' and this causes serious difficulties.

\section{The finite-dimensional approximation of $P^+$}

We already mentioned in the footnote on p. \pageref{space} that we
have to deal with good models for the spaces $P^+$ and $\Omega^+$,
and now we recall one of the approaches which was used
by Morse \cite{Morse} (see also
\cite{Milnor}).

Let us sketch his constructions and results.

Let $M$ be a compact Riemannian manifold. Let $\delta$ be a positive
constant such that any two points $p, q \in M$ with $d(p,q) \leq
\delta$ \footnote{Here $d(p,q)$ is the distance function
corresponding to the Riemannian metric on $M$.} are connected by a
unique geodesic $\gamma(p,q)$ with length $\leq \delta$ and moreover
this geodesic smoothly depends on $p$ and $q$. For $M$ a compact
manifold such a constant $\delta$ always exists and for it we can
take any positive constant which is less than the injectivity
radius.

Let $P_a^+$ be a space of closed curves with marked points and
length less than $a$ and let
$$
N = \left[\frac{a}{\delta}\right] + 1.
$$
Starting from the marked point $r(0)$ divide every curve into $N$
successive pieces of equal length. If we take the length parameter
on $r(t)$, then the separating points are
$$
r(0), r(1/N), \dots, r((N-1)/N), r(1)=r(0),
$$
Then let us replace every curve segment between $r(i/N)$ and
$r((i+1)/N)$, $i=,\dots,N-1$, by a unique shortest geodesic which
connects these points, i.e. by $\gamma(r(i/N),r((i+1)/N))$. By
construction, its length is less than $\delta$. Then the curve $r$
is replaced by a piecewise geodesic curve. All such geodesic
polygons form a subset
$$
\Pi_N \subset \underbrace{M \times \cdot \times M}_{N}
$$
which is the closure of an open domain $U$ in the $N$-th power of
$M$. The length functional gives rise to a smooth function $f$ on
$U$, critical points of $f$ lie in $U$ and they are exactly closed
geodesics. We have
$$
\dim U = Nd, \ \ \ \mbox{where $d = \dim M$}.
$$
The points in $\Pi_N$ are parameterized by sequences of the
end-points by geodesic segments:
$$
(x_0,\dots,x_{N-1}), \ \ \ x_i \in M.
$$
In these terms the action of $\sigma$, the reversion of
an orientation, reduces to
$$
(x_0,\ldots,x_{N-1}) \stackrel{\sigma}{\longrightarrow}
(x_{N-1},\ldots,x_0).
$$
The group $S^1$ contains the subgroup
$\Z_N$ which acts as follows
$$
(x_0,\dots,x_{N-1}) \stackrel{\theta}{\longrightarrow}
(x_\theta,\dots,x_{(N-1+\theta)\mathrm{mod}\,N}), \ \ \ \theta \in
\Z_N.
$$

The function $f$ is invariant with respect both to the
$\Z_N$ action and to $\sigma$.

Even after passing to $\Pi_N/\Z_N$
any nontrivial closed geodesic corresponds to an
$S^1$-family of critical points of $f$.

However now we are in a finite-dimensional situation and can define
the index and the nullity of a closed geodesic as of a critical
point of smooth function. The index a closed geodesic was defined by
Morse as follows. Let $x=(x_0,\ldots,x_{N-1})$ corresponds to a
closed geodesic $g$. At every point $x_i$ let us take an
$(d-1)$-dimensional (i.e. a hypersurface) open submanifold $W_i
\subset M$ which passes through $x_i$ and transversal to the
geodesic $g$. Assuming that $W_i$ are sufficiently small, we have an
embedding
\begin{equation}
\label{w}
W=W_0 \times \ldots \times W_{N-1} \subset \Pi_N \subset
M^N.
\end{equation}

The restriction of the length function onto $W$ is a smooth
function
$$
f: W \to \R
$$
and we define the index and the nullity of the
critical point $x$ as in a finite-dimensional case:

\begin{itemize}
\item
{\it the index} of $x$ equals the index of the Hessian of $f$ at $x$,
i.e. the dimension of a tangent subspace at $x$ on which the Hessian of
is negatively-defined;

\item
{\it the nullity} of $x$ is the difference of the dimensions of the
tangent space in $x$ and its maximal subspace onto which the
restriction of the Hessian is nondegenerate;

\item
a geodesic is {\it nondegenerate} if its
nullity vanishes.
\end{itemize}

There is another way to define nondegenerate closed geodesics and
their indices proposed by Bott \cite{Bott54}. Given a smooth
function $f:X \to \R$ on a finite-dimensional manifold and a compact
connected submanifold $Y\subset X$ which is a critical set of $f$,
we say that
$Y$ is nondegenerate if the nullity of $f$ at every point $x\in Y$
equals to the dimension of $Y$. The index of $Y$ is defined as the
index of $f$ at $x \in Y$. Now it is said that a geodesic $g=x \in
\Pi_N$ is nondegenerate if the corresponding $S^1$-family of
critical points is nondegenerate in the Bott sense and the index of this
geodesics is the index of the corresponding critical submanifold homeomorphic to
$S^1$.
\footnote{A function for which all critical sets are nondegenerate is called
a {\it Morse--Bott function} due to the Bott article
\cite{Bott54} from 1954.

We have to stress that it was Pontryagin
\cite{Pontryagin,Pontryagin39} who at the first time considered such
functions and even applied the extended Morse theory for solving a
topological problem: a computation of the Betti numbers of classical
Lie groups in 1935. In this case on every such a group there is a
function with only two critical levels. For instance, given the
group $SO(N)$ of orthogonal $N\times N$-matrices Pontryagin
considered the entry $f=a_{11}$ in the top left corner of an
orthogonal matrix as a smooth function on $SO(N)$ and proved that
its only critical levels are $\pm 1$ and the corresponding critical
sets are naturally diffeomorphic to $SO(N-1)$. Carefully analyzing
how normal neighborhoods of critical sets are glued by the gradient
flow of $f$ and starting with the simplest case $SO(3)= \R P^3$,
Pontryagin computed successively the $\mod \, 2$ Betti numbers of
all the groups $SO(N)$ and did the same for three other classical
series.

However Pontryagin did not consider the general situation and did not define
nondegenerate critical submanifolds and their indices.}

\section{The type numbers of a closed geodesic in $\Omega^+$}

The sets of type (\ref{w}) give only local finite-dimensional
approximations of the space $\Omega^+ = P^+/S^1$. The idea of Morse
which is carefully discussed on p. 249 is to use the limit of
$\Pi_N$ as $N \to \infty$ for describing $\Omega^+$. If we have a
geodesic polygon $x \in \Pi_N$ we may split every edge of it into
$l$ parts of equal length. Thus we obtain a polygon with $Nl$ edges
and this canonical procedure gives us projections:
$$
\Pi_N \to \Pi_{Nl} \ \ \ \mbox{for all $N,l$.}
$$
To these projections there corresponds a directed system of
relative homology groups $H_\ast(\Pi_N,M)$ where by $M$ we denote the
subset formed by trivial (one-point) curves.
{\it The $k$-th circular connectivity}
as it is defined by Morse \cite{Morse}
is the rank of the group
$$
\bar{p}_k = \rank \lim_\to H_k(\Pi_N/\sigma,M\;\Z_2),
$$
and {\it the $k$-th sensed
circular connectivity}
$$
p_k = \rank \lim_\to H_k(\Pi_N;\Z_2)
$$
was introduced by Bott  \cite{Bott54}.
Since, for any group of coefficients we have
$$
\lim_\to H_i(\Pi_N,M) = H_i(\Omega^+,M),\ \ \ \ \
\lim_\to H_i(\Pi_N/\sigma,M) = H_i(\Omega,M)
$$
(see Theorem 9 in \cite{Schwarz}),
these topological invariants give lower estimates for sums of the
type numbers of critical sets of the length functional, and hence
imply the existence of nontrivial closed geodesics.

Morse stated that
{\sl the type numbers $m_k$ (over $\Q$ or $\Z_p$) of a nondegenerate
closed geodesic of index $i$ in $\Omega^+$ (or $\Omega$)
are equal to}
\begin{equation}
\label{pitype}
m_k = \delta_k^i
\end{equation}
\cite[p. 296]{Morse}. As it was first mentioned in \cite{Schwarz} in
general this is false however this is valid for prime closed
geodesics. The failure of (\ref{pitype}) is explained by the fact
that

\begin{itemize}
\item
{\sl the $SO(2)$ action is not free at iterated closed curves.}
\end{itemize}

The right answer is given by the following
\footnote{The final corrected version of this result
was given in \cite{Alber}.}

\begin{theorem}[\cite{Schwarz,Alber}]
\label{th1}
Let $g$ be a prime closed geodesic, let $h=g^s$ be its $s$-th
iteration, let $i(d)$ be the index of $g^d$, and let $d_p$ be the
greatest divisor of $s$ which is not divided by a prime number $p$.
We assume that the geodesic $g$ and its iterates are nondegenerate.

Then the type numbers $m_k$ and $m_k^{(p)}$ of $h=g^s$ in $\Omega^+$ (or
in $\Omega = \Omega^+/\sigma$) with respect to coefficients in
rational numbers and in $\Z_p$ take the form

\begin{enumerate}
\item
if $s=1$, then
$$
m_k = m_k{(p)} = \delta_k^{i(1)};
$$

\item
if $s \geq 2$, then

\begin{enumerate}[(a)]
\item
if $s$ is odd or both $s$ and $i(2)-i(1)$ are even, then
$$
m_k = \delta_k^{i(s)},
$$
$$
m_k^{(p)} =
\begin{cases}
1 & \text{if $i(d_p) + 2 \leq k \leq i(s)$}, \\ 0 & \text{otherwise};
\end{cases}
$$

\item
if $s$ is even and $i(2)-i(1)$ is odd, then
$$
m_k = m^{(p)}_k = 0 \ \ \ \mbox{for all $k$ and $p \neq 2$},
$$
$$
m_k^{(2)} =
\begin{cases}
1 & \text{if $i(d_2) +2 \leq k \leq i(s)$}, \\ 0 & \text{otherwise}.
\end{cases}
$$
\end{enumerate}
\end{enumerate}
\end{theorem}

The formula (\ref{pitype}) was used in \cite{Morse,Bott54} for computing
$\mod 2$\, homologies of the spaces $\Omega$ and $\Omega^+$, respectively,
for spheres. The failure of this fact explains why these computations were
not correct and, to our knowledge, these homologies are
still not known.
\footnote{See the discussion in \cite{Hingston}.}

The rational cohomologies of these spaces were computed in \cite{Schwarz}.
We formulate these results in terms of the Poincare series
$$
p(X,t) =  \sum_{k \geq 0} \rank H^k(X;\Q) \cdot t^k, \ \
p(X,Y,t) = \sum_{k \geq 0} \rank H^k(X,Y;\Q) \cdot t^k.
$$

\begin{theorem}[\cite{Schwarz}]
\label{th2}
For the spaces $\Omega^+(S^n)$ and $\Omega(S^n)$, we have

\begin{enumerate}
\item
for $n$ even
$$
p(\Omega^+,M,t) = t^{n-1}\left(
\frac{1}{1-t^2} + \frac{t^{2n-2}}{1-t^{2n-2}}\right),
$$
$$
p(\Omega,M,t) = t^{n+1}\left(
\frac{1}{1-t^4} + \frac{t^{2n-4}}{1-t^{4n-4}}\right);
$$

\item
for $n \geq 3$ odd
$$
p(\Omega^+,M,t) = t^{n-1}\left(
\frac{1}{1-t^2} + \frac{t^{n-1}}{1-t^{n-1}}\right),
$$
$$
p(\Omega,M,t) = t^{n+1}\left(
\frac{1}{1-t^4} + \frac{t^{n-3}}{1-t^{2n-2}}\right);
$$

\item
$$
p(\Omega^+,t) = 1 - t^{n+1} + p(\Omega^+,M,t),
$$
$$
p(\Omega,t) = 1 - t^{n+1} + p(\Omega,M,t)
$$
for all $n$.
\end{enumerate}
\end{theorem}

For this computations Schwarz used the Leray--Serre spectral sequence
applied to the homotopy quotients
$$
X_{S^1} = X \times_{S^1} ES^1, \ \ X = P^+,
$$
and their (homotopy) fibrations
$$
X_G \stackrel{X}{\longrightarrow} BG, \ \ \ G=S^1.
$$
Since
$$
H^\ast(\Omega^+;\Q) = H^\ast_{S^1}(P^+;\Q),
$$
that enables us to compute the rational cohomologies of
$\Omega^+=P^+/S^1$. This method was just introduced by Borel
\cite{Borel} and now
$$
H^\ast(X_G) = H^\ast_G(X)
$$
are called the equivariant cohomology of a space $X$ with $G$-action.
Later $H^\ast(\Omega^+(M),M;\Q)$ were computed by the same method for
other compact rank one symmetric spaces $\C P^n, \H P^n$, and $\mathrm{Ca} P^2$
by Hingston \cite{Hingston} and, finally, Rademacher \cite{Rademacher} did that
for all spaces $M$ for which $H^\ast(M;\Q)$ is a truncated polynomial ring
in one variable.

\section{The Bott formula for the Morse indices of iterated closed geodesics
and its applications}

Before demonstrating an important fact, i.e. Theorem \ref{th4},
which was proved by using Theorems \ref{th1} and \ref{th2}, we
expose the Bott formula \cite{Bott56}. It plays a crucial role in
proving the existence of a few closed geodesics and, especially, in
the known proofs of the existence of infinitely many (geometrically
distinct) closed geodesics on certain manifolds (see, in particular,
Theorem \ref{th5}). It is also necessary for proving Theorem
\ref{th4}.

\begin{theorem}[\cite{Bott56}]
\label{th3}
Let $g$ be a prime closed geodesic on a $d$-dimensional manifold $M$.
Then there non-negative integer-valued functions $\Lambda_g$ and $N_g$
on the unit circle $|z|=1$ such that

\begin{enumerate}
\item
$N_g(z) = \dim \mathrm{ker}\, (P_g-z)$
where $P_g$ is the Poincare map of the geodesic $g$ (it is given by
$(d-1)\times (d-1)$-matrix);

\item
$\Lambda_g(z)$ is constant at points at which $N_g(z)=0$ and its
jumps are always bounded in absolute value by $N_g(z)$;

\item
$\Lambda_g(z) = \Lambda_g(\bar{z}), N_g(z) = N_g(\bar{z})$;

\item
the indices $i(g^n)$ and the nullities $\nu(g^n)$
of iterated geodesics $g^n$ are equal to
$$
i(g^n) = \sum \Lambda_g(\omega),  \ \ \
\nu(g^n) = \sum N_g(\omega),
$$
where $\omega$ ranges over the $n$-th roots of $1$ if the parallel
translation along $g$ does not change the orientation of $M$ and
over the $n$-th roots of $-1$ overwise.
\end{enumerate}
\end{theorem}

These formulas for indices and nullities are particular cases
of facts concerning general self-adjoint differential equations with
periodic coefficients.

A delicate corollary of Theorems \ref{th1} and \ref{th3}, which is
useful in applications, is as follows.

\begin{lemma}
[\cite{Fet}]
\label{lemma1}
Let $g$ be a prime orientation-preserving
closed geodesics of index $i(g)=l$ and let $g$ and all its iterations
$g^n, n=1,2,\dots$, are nondegenerate.
Then for all $n$ we have
$$
m_k(g^n) = 0 \ \ \ \ \mbox{if $k-l$ is an odd integer},
$$
where $m_k$ are the type numbers in $\Omega^+$ with respect
to rational coefficients.

Hence, if $m_k(g^l) \neq 0$, then $k-i(g)$ is even.
\end{lemma}

{\sc Proof.} By Theorem \ref{th1}
$$
m_k(g^n) =
\begin{cases}
\delta_k^{i(g^n)} & \text{when $i(g^2)-i(g)$ is even} \\
0 & \text{otherwise}.
\end{cases}
$$
However, by Theorem \ref{th3}, if $i(g^2)-i(g)$ is even, then
for all $n$ we have $i(g^n)-i(g)$ is even and Lemma follows.

From Lemma \ref{lemma1}
Fet derived the following theorem whose proof substantially
uses Theorem \ref{th2}.
\footnote{In \cite{Fet} it is claimed the existence of the second closed
geodesic on a non-simply-connected manifold
is also proved but the arguments are wrong and this question
is still open. This is explained on p. \pageref{fetmistake}.}

\begin{theorem}[\cite{Fet}]
\label{th4} Let $M$ be a simply-connected closed Riemannian
manifold. Assume that all closed geodesics are nondegenerate. Then
there are closed geodesics $g_1$ and $g_2$ with $m_k (g_1) \neq 0$
and $m_{k+1}(g_2) \neq 0$ with some $k < \dim M$. In particular,
these geodesics are geometrically distinct.
\end{theorem}

The proof of Theorem \ref{th4} is based on the fact that for spheres
we have
$$
H_i(\Omega(S^n),S^n;\Q) = 0 \ \ \ \mbox{for all $i < n$}
$$
which, in its turn, is not true for the spaces of oriented curves:
$$
\rank H_{n-1}(\Omega^+(S^n),S^n;\Q) = 1.
$$

One of the most substantial results which is derived by using the
Bott formulas is the following theorem of Gromoll and Meyer:

\begin{theorem}
[\cite{GM}]
\label{th5}
Let $M$ be a simply-connected closed Riemannian manifold.
If the Betti numbers of the free loop space $\Lambda M$ (for any given field
of coefficients) are unbounded then there
are infinitely many geometrically distinct (nontrivial) closed geodesics.
\end{theorem}

Since the Morse theory establishes a strong relation between
critical points and the Betti numbers, for a proof Gromoll and Meyer
did show that finitely many closed geodesics and their iterations
can not generate the homology of the free loop space with unbounded
Betti numbers.

More delicate consequences of the Bott formula were found and used
by Rademacher (the average index formula \cite{Rademacher}) and Long
(the common index jump theorem \cite{LongZhu}).

\section{On lower estimates for the number of closed geodesics}

For completeness we expose here the main results concerning the existence
of closed geodesics obtained by using variational methods and, in particular,
by methods and ideas introduced by Morse \cite{Morse}.

We remark that whenever in this section we speak about
estimates for the number of closed geodesics we mean prime closed geodesics.

\subsection{Geodesics on non-simply-connected manifolds}

For a nonsimply-connected manifold $M$ to every element $\tau \in
[S^1,M]$, the set of free homotopy classes of maps from $S^1$ into
$M$, there corresponds a connected component $P^+_\tau(M)$. Let
$c_\tau$ be the infimum of the lengths of curves from $P^+_\tau(M)$
for some nontrivial class $\tau$. Then consider the
finite-dimensional approximation $\Pi_{N,\tau}$ of the subset
$\{\mathrm{length} < c +\varepsilon\} \in P^+_\tau$, where
$\varepsilon$ is any positive constant. On this finite-dimensional
set the length reduces to a function which is smooth in the interior
of $\Pi_{N,\tau}$. It is easy to show that this function admits its
minimum at a smooth point which, by construction, is a closed
geodesic $\gamma_\tau$ with minimal length in $\tau \in [S^1,M]$.
Since $\tau$ is nontrivial, the geodesic $\gamma_\tau$ is also
nontivial. This proves

\begin{theorem}
[Hilbert]
\label{th6}
On every non-simply-connected Riemannian closed manifold there exists a closed
geodesic which is not homotopic to zero.
\end{theorem}

Now we come to the problem: how to prove that the geodesics which are minimal
in the classes $\tau_1$ and $\tau_2$ are geometrically different?

We recall that
$$
[S^1,M] = \pi_1/\sim, \ \ \ \mbox{where $x\sim y$ if there is $z$ such that
$x=zyz^{-1}$}
$$
and $H_1(M)$ is the abelianization of $\pi_1$:
$$
H_1 = \pi_1/[\pi_1,\pi_1].
$$
Let
$$
H_1(M) = \Z + \Z + A
$$
where $A$ is some group and let $a$ and $b$
be the generators of the first two $\Z$ subgroups. Let us take the free
homotopy classes $\tau_p$, where $p$ range over prime numbers,
which represent the homology classes $a + pb$. It is easy to notice that
the closed geodesics $\gamma_{\tau_p}$ are geometrically distinct.
In particular, this argument shows that metrics on all closed oriented
surfaces except $S^2$ have infinitely many closed geodesics.
\footnote{By passing to the universal covering, we may claim the same for
all non-oriented surfaces except $\R P^2$.
The cases of $S^2$ and $\R P^2$ are covered by Theorem \ref{th9}.}

Until recently for a general fundamental group (save of cyclic
groups) it is not proved that the minimizers in different free
homotopy class give us at least two geometrically distinct closed
geodesics! In fact, we can not exclude the existence of a
finitely-presented group $\pi$ with the following property:

\begin{itemize}
\item
there is an element $a \in \pi$ such that
every nontrivial element $b \in \pi$
is conjugate to some power of $a$: $a^k \sim b$.
\end{itemize}

At least it is known that such noncyclic finitely-generated groups
do exist \cite{Guba} (however all known such examples are not
finitely-determined groups). This demonstrates a flaw in the
arguments from \cite{Fet}. Despite the claim from this article the
existence of at least two closed geodesics on every
non-simply-connected closed manifold is not proved yet.
\label{fetmistake}

For seeking more closed geodesics we have to look for non-minimal
closed geodesics. For proving their existence one can use higher
homotopy groups of $\Omega^+_\tau(M)$ computed in
\cite{Ballmann,Taimanov85} in terms of action of the fundamental
group on higher homotopy groups of $M$. Until recently the class of
finitely-presented groups is not not well understood and that
obstructs deriving very general statements on the existence of many
closed geodesics on non-simply-connected manifolds. \footnote{ We
would like to mention another approach to study the existence of
closed geodesics in which algorithmic properties of the fundamental
group plays the main role. In 1976 Gromov \cite{Gromov} stated the
following theorem whose proof was exposed in \cite{Nabutovsky}:

\begin{itemize}
\item
{\sl if the word problem for $\pi_1(M)$ is unsolvable, then on a
closed manifold $M$ there exists an infinite family of prime
contractible closed geodesics $g_n$ such that each of them provides
local minimum for the length functional and $\mathrm{length}\, g_n
\to \infty$ as $n \to \infty$.}
\end{itemize}

The proof shows that if there are only finitely many
prime contractible closed geodesics then there is an algorithm to
decide whether any word in $\pi_1$ is trivial or not.

Nabutovsky distinguished in terms of the Kolmogorov
complexity of the word problem the class $K$ of groups such that for a closed
manifold $M$ with $\pi_1(M) \in K$ there exists a constant $c=c(M) > 1$
such that the number $N(t)$ of contractible closed geodesics
of length $\leq t$ grows
at least as $c^t$ \cite{Nabutovsky}.}

However for $\pi_1=\Z$ Bangert and Hingston
by very nice arguments involving the actions of $\pi_1$ on $\pi_n, n\geq 2$,
proved not only the existence of
infinitely many closed geodesics but also the asymptotic formula for the growth
of their lengths \cite{BH}. Let $N(t)$ be a number of geometrically distinct
closed geodesics of length $\leq t$. Then
$$
N(t) > C \frac{t}{\log t}
$$
with $C$ a positive constant.

\subsection{Geodesics on simply-connected manifolds}

{\sc On the existence of a closed geodesic.}

The existence of a nontrivial closed geodesic on a sphere $S^n$
of any dimension $\geq 2$ was established by Birkhoff as follows.

Let us take a singular foliation of the sphere by circles and
points. For a two-dimensional sphere $S^2$ realized in $\R^3$ as the
unit sphere which foliation is formed by plane sections where planes
are orthogonal to some fixed axis. Assume that these curves are
oriented, then this family of closed curves on $S^n$ form an
$(n-1)$-dimensional cycle $\beta_{n-1} \in P^+(S^n)$. Since the
family is glued together into a manifold which is mapped onto $S^n$
with degree one, this cycle is nontrivial in $H_{n-1}(P^+(S^n);\Q)$.
It generates also a cycle in $H_{n-1}(P(S^n);\Z)$. Since in
non-oriented case the degree of the map is defined $\mod 2$ this
cycle may lie in $2$-torsion in $H_{n-1}(P(S^n);\Z)$ and, in fact,
that holds:
$$
2 \beta_{n-1} = 0 \in H_{n-1}(P(S^n);\Z).
$$
However from (\ref{morseineq}) it follows that
$\beta_{n-1}$ suspends on some closed geodesic.

The existence of a nontrivial closed geodesic
on any simply-connected closed manifold is derived by using these
Birkhoff cycles. By the Hurewicz Theorem, for any simply-connected
closed manifold
$M$ there are nontrivial higher homotopy groups. In particular, if
$$
\pi_1(M) = 0, \ \ \ H_i(M) = 0 \ \mbox{for $i<k$ and} \ \ H_k(M) \neq 0,
$$
then
$$
\pi_i(M) = 0 \ \ \mbox{for $i<k$ and} \ \ \pi_k(M) = H_k(M).
$$
Let us take a mapping $f: S^k \to M$ which realizes a nontrivial element
$[f] \neq 0$ in $\pi_k(M)$. We take a singular foliation of $S^k$ and
mapping it into $M$ we obtain the image $f_\ast(\beta_{k-1}) \in
P^+(M)$ of the Birkhoff cycle. This cycle is also nontrivial,
for the same reasons as the Birkhoff cycle, and on this cycle there
suspends a closed geodesic.

Together with Theorem \ref{th6} this implies

\begin{theorem}
[\cite{LF}]
\label{th7}
On every closed Riemannian manifold there exists a clo\-sed geodesic.
\end{theorem}

\smallskip

{\sc On the existence of infinitely many closed geodesics.}

Now there is a strong belief that on every, at least simply-connected, closed
Riemannian manifold there are infinitely many (geometrically distinct)
closed geodesics. A substantial progress in proving that was achieved in
1960-90s. We shall expose the main known results.

These results are proved mainly by applying the Morse inequalities
(\ref{morseineq}) for establishing that the Betti numbers
of $P(M)$ or $\Omega(M)$ can not be generated by finitely many
nontrivial closed geodesics and their iterates. More complicated
ring structure of the cohomologies of these spaces is not used.
Other approaches were used only for two-dimensional manifolds
(see the survey \cite{Taimanov}).

The development of rational homotopy theory (minimal models by Sullivan)
led to the following result by Sullivan and  Vigu\'e-Poirret:

\begin{theorem}
[\cite{Sullivan}]
\label{th8}
If the rational cohomology ring of a closed simp\-ly-con\-nected manifold
$M$ is not a ring of truncated polynomials of one variable, then
the rational Betti numbers of the free loop space $\Lambda M$ are unbounded.
\end{theorem}

Together with Theorem \ref{th5} this implies

\begin{corollary}
[\cite{Sullivan}]
\label{cor1}
If the rational cohomology ring of a closed simp\-ly-con\-nected manifold
$M$ has at least two multiplicative generators, then any Riemannian metric
on $M$ has infinitely many closed geodesics.
\end{corollary}

So we left only with manifolds for which
\begin{equation}
\label{truncated}
H^\ast(M;\Q) = \Q[u]/u^{d+1},
\end{equation}
which are, for instance, compact rank one symmetric spaces.
For them the rational Betti numbers of $\Lambda M$ are bounded.

In early 1990s its was proved that

\begin{theorem}
[\cite{Bangert,Franks}]
\label{th9}
Every Riemannian metric on a two-sphere has
infinitely many closed geodesics.
\end{theorem}

Theorem \ref{th9} is proved by combining results of Bangert
\cite{Bangert} and by Franks \cite{Franks}. It needs to consider two
cases: 1) there exists a simple closed geodesic for which the
Birkhoff map is defined; 2) there is no such a geodesic. In the
first case Franks proved the statement by using the theory of
area-preserving diffeomorphisms which goes back to the
Poincare--Birkhoff theorem, and in the second case Bangert showed
that there exists a simple closed geodesic without conjugate points
(when considered as geodesic defined for all times) and then used
the Morse theory to derive from that the existence of infinitely
many closed geodesics. In \cite{Hingston93} the arguments by Franks
were replaced by Morse-theoretical reasonings however it is still
necessary for proof to consider two cases independently and the
difference between them is formulated in dynamical terms.

For other manifolds the only known lower estimate which differs from
just the existence of a closed geodesic was
recently found for three-spheres:

\begin{theorem}
[\cite{Long}]
\label{th10}
Every Riemannian metric on a three-sphere has at least two closed geodesics.
\end{theorem}

\smallskip

{\sc On lower estimates for the number of closed geodesics of
ge\-ne\-ric metric.}

It was Morse who first started to consider metrics meeting some
generic conditions. He considered metrics for which all closed
geodesics are nondegenerate (these metrics are called now {\it
bumpy} metrics). In this case for every closed geodesic only one of
its type numbers (with respect rational coefficients) does not
vanish. This property is $C^r$-generic for $2 \leq r \leq \infty$:
that was announced by Abraham in 1960s and the proof of that was
given by Anosov in \cite{Anosov82}.

Theorem \ref{th4} gave the best known lower estimate for the number of
closed geodesic of a bumpy metric which is valid for all
simply-connected closed manifolds.

However one may consider other generic conditions.

The result by Klingenberg and Takens \cite{KT} reads that

\begin{itemize}
\item
the condition
\footnote{We shall call it the KT condition.}
for a metric that
either all closed geodesics are hyperbolic,
\footnote{In particular,
they are nondegenerate and
for them index grows linearly: $i(g^n)=ni(g)$.} either
there exists a nonhyperbolic closed geodesic of twist type
\footnote{In this case,
by Moser's result, near this geodesic there exists infinitely many
closed geodesics. The proof of that
relates to the same theory as the
result by Franks used in proving Theorem \ref{th9}.}
is a $C^r$-generic condition for $4 \leq r \leq \infty$.
\end{itemize}

Until recently there are no known examples of simply-connected
closed Riemannian manifolds with only hyperbolic closed geodesics.
However to prove the existence of infinitely many closed geodesics
for metrics meeting the KT condition it is enough to consider
manifolds with only hyperbolic geodesics.

For non-simply-connected manifold that was done by Ballmann, Thorbergsson,
and Ziller who proved that

\begin{theorem}
[\cite{BTZ}] \label{th11} Given a closed manifold $M$, if there is a
nontrivial element $a \in \pi_1(M)$ such that $[a^m]=[a^n] \in
[S^1,M]$ for some integers $m \neq n$, then for any metric on $M$
meeting the KT condition there is a constant $C>0$ such that
$$
N(t) > C \frac{t}{\log t}
$$
where $N(t)$ is the number of geometrically distinct closed
geodesics of length $\leq t$.
\end{theorem}

The hypothesis of Theorem \ref{th11} is satisfied if
$a \neq 1$ and $a$ is of finite order, i.e. for groups with torsion elements.

The proof of the existence of infinitely many closed geodesics
in the situation of Theorem \ref{th11} is as follows. Let us assume that
all geodesics are hyperbolic and  let us
take a minimal geodesic $g$ in class $[a]$. Then for all $l$ the geodesics
$g^{lm}$ and $g^{ln}$ have index $0$ and are freely homotopic to each other.
Since they lie in a connected space to every such a pair $g^{lm}$ and
$g^{ln}$ there corresponds a geodesic of index $1$ which suspends on
a one-cycle which joins $g^{lm}$ and $g^{ln}$. These geodesics of index $1$
are prime since for hyperbolic geodesics $i(h^n) = ni(h)$.

The case of KT metrics on simply-connected manifolds was considered
by Rademacher who computed the rational
cohomologies of $\Omega^+(M)$ for all manifolds meeting
(\ref{truncated}) and
showed that they can not be generated by finitely many hyperbolic geodesics
and their iterates. For proving that he introduced and used the average
index formula for closed geodesics.

\begin{theorem}
[\cite{Rademacher}]
\label{th12}
Every metric on a simply-connected closed manifold $M$ which meets
the KT condition has infinitely
many closed geodesics.
\end{theorem}

Later Rademacher considered so-called {\it strongly bumpy metrics}
which are defined as metrics for which all (non-trivial) closed
geodesics are nondegenerate, the eigenvalues of the Poincare maps
for them are simple and every finite set of the eigenvalues of norm
one are algebraically independent. In \cite{Rademacher94} he proved
that this condition is $C^r$-generic for $2 \leq r \leq \infty$ and
by using the average index formula also proved that for such metrics
on simply-connected manifolds there exist infinitely many closed
geodesics.

\section{Finsler and magnetic geodesics}

Since our task was mostly to complement the new edition of Morse's book
by necessary comments, we did not write an extensive
survey on applications of
the Morse theory to closed geodesics. We also put off
the theory by Lusternik and Schnirelmann
who used different topological invariants (category, cup length)
which now serve as very strong tools in nonlinear analysis and calculus of
variations. However the proof of their main theorem was completed only recently
\cite{Ballmann78} (see also \cite{Taimanov}).

At the end we would like to mention briefly two active subjects where
the Morse theory is being applied and where there are still a plenty of
open problems.

\subsection{Finsler geodesics}

Although the study of closed geodesics for Finsler metrics was initiated
by Anosov in the middle of 1970s \cite{AnosovICM} it became an active area
of research only recently.

A Finsler metric $F(x,\dot{x})$
is a nonnegative function on the tangent bundle which is convex in
$\dot{x}$ and meets the condition
$$
F(x,\lambda \dot{x}) = \lambda F(x,\dot{x}) \ \ \ \mbox{for $\lambda > 0$}.
$$
Riemannian metrics supply a particular class of Finsler metrics:
$$
F(x,\dot{x}) = \sqrt{g_{ik}(x) \dot{x}^i\dot{x}^k} = |\dot{x}|.
$$
In fact, the whole Morse theory works for Finsler metrics without
changes: it works for sufficiently general functionals on loop
spaces and, in particular, to the Finsler length. However we have to
distinguish between reversible and irreversible Finsler metrics: a
Finsler metric is called reversible if
$$
F(x,\dot{x}) = F(x,-\dot{x}).
$$

If a metric is reversible, then the change of orientation of a
closed curve preserves the Finsler length and we may deal with the
space of non-oriented curves, i.e. with $P(M)$ and $\Omega(M)$.
Otherwise, for irreversible metrics we have to look for critical
points in $P^+(M)$ and $\Omega^+(M)$.

There are strong differences between the homologies of spaces of
oriented and of spaces of un-oriented curves. For instance, the
Birkhoff cycle $\beta_n$ is nontrivial in $H_n(P^+(S^{n+1});\Q)$ and
vanishes in $H_n(P(S^{n+1});\Q)$. For instance, this explains why
Theorem \ref{th4} does not hold for geodesics of irreversible
metrics (the proof from \cite{Fet} works only for reversible metrics
and for irreversible metrics the counterexample is given by the
Katok example of a Finsler metric on $S^2$ for which there are only
two closed geodesics $g_1$ and $g_2$ with $m_k(g_1) \neq 0$ and
$m_{k+1}(g_2)\neq 0$ \cite{Katok}). However, the proofs of Theorems
\ref{th5}, \ref{th6}, \ref{th7}, \ref{th11}, \ref{th12} and of
Corollary \ref{cor1} work without changes for all Finsler metrics.

In general we can not expect the existence of infinitely many closed
Finsler geodesics. An example of Katok which we already mentioned
shows that there is an irreversible bumpy Finsler metric on $S^2$
admitting only two prime closed geodesics which are given by the
same curve passed in different directions \cite{Katok,Ziller}.

In the last years the periodic problem for Finsler geodesics
as well the differential geometry in large of Finsler spaces
attracts a lot of attention \cite{BangertLong,LongEMS,Rademacher04}.
In particular, very recently Bangert and Long proved
the existence of at least
two geometrically distinct closed geodesics on every Finsler $2$-sphere
\cite{BangertLong} (this was claimed in \cite{AnosovICM} but
the proof did never appear).

An excellent survey on the periodic problem in the Finsler geometry
and a list of open problems is presented in
\cite{LongEMS}.

\subsection{Magnetic geodesics}

A simple example of a Finsler metric is provided
by the following function
$$
L(x,\dot{x}) = |\dot{x}| + A_i \dot{x}^i,
$$
where $A_i dx^i$ is a one-form.
If $A$ is sufficiently small:
$$
L(x,\dot{x}) > 0 \ \ \ \mbox{for $\dot{x} \neq 0$},
$$
we have an irreversible Finsler metric which in the Finsler geometry
is called a Randers metric and we may apply the classical Morse theory.

We note that $A$ comes into the corresponding
Euler--Lag\-ran\-ge equations via its differential
$$
F = \sum_{i<k} F_{ik} dx^i \wedge dx^k, \ \ \
F_{ik} = \frac{\partial A_k}{\partial x^i} -
\frac{\partial A_i}{\partial x^k},
$$
and that is easily explained by physics where
solutions of the Euler--Lagrange equations
for the Lagrangian
$$
L(x,\dot{x}) = \sqrt{E g_{ik}\dot{x}^i\dot{x}^k} + A_i \dot{x}^i
$$
are (up to a parameterization)
trajectories of a charged particle in a magnetic field
on a fixed energy level $E$, where the Riemannian metric determines the
kinetic energy, and
a closed two-form $F$ is a magnetic field \cite{Novikov}.
Hence for physical reasons only $F$
has to be globally defined and if $F$ is not an exact form, then the one-form
$A$ is defined only locally.
So we have two cases which are drastically different from periodic
problems considered before by the Morse theory:

a) the form $A$ is defined only locally, in this case the formula
\begin{equation}
\label{magnetic}
S(\gamma) = \int_\gamma (|\dot{x}| + A_i^\alpha \dot{x}^i) dt
\end{equation}
determines a multi-valued functional on the space of closed curves
(for instance, on $P^+(M)$) with a uniquely and globally defined
variational derivative $\delta S$ (here $A_i^\alpha dx^i$ is the
form $d^{-1} F$ defined in some domain $U_\alpha$ which contains
$\gamma$ and for which the second Betti number vanishes);

b) even if $A$ is globally defined, then the functional
(\ref{magnetic}) may be not bounded from below. In this situation when using
the Morse theory
we have to consider the homology groups
$H_\ast(P^+(M),\{S \leq 0\})$ instead of $H_\ast(P^+(M),M)$ in
the classical Morse theory.

The extremals of the functional
(\ref{magnetic}) are called {\it magnetic geodesics}.
Sometimes the flow of such extremals is called {\it twisted geodesic flow}
however we do not use this term and
prefer to stress the physical origin of the subject.

To the functional (\ref{magnetic}) one may apply the whole local
part of the Morse theory,
defining the index and the nullity of extremals.

The study of the existence of critical points of these functionals
as well as developing the analog of the Morse theory for one-forms
(multivalued functions with globally defined differential) was
started by Novikov in 1980
\cite{Novikov81,Novikov-a,Novikov-b,Novikov}. The finite-dimensional
theory was far-developed \cite{Farber,Pazhitnov} and its main
counterparts such as the Morse--Novikov inequalities and the Novikov
ring found many applications. The ideas concerning multi-valued
functions found applications in mathematical physics. Speaking about
closed extremals in magnetic fields, we note that topological
problems concerning the lower estimates for the number of critical
levels were also basically solved (Novikov's principle of ``throwing
out cycles'') \cite{Novikov,Taimanov83,Taimanov84}. However the main
difficulty appeared in the analytical part.

For deriving the existence of critical points
from the Morse inequalities (\ref{morseineq})
it needs to have a good deformation which
decreases the value of $S$ at noncritical points
and lets a nontrivial relative cycle to suspend on
a critical point of $S$ in the limit. For the length functional Morse uses
the following deformation on $\Pi_N$: given a geodesic $N$-gon
$\gamma = (x_0,\dots,x_{N-1}) \in \Pi_N$,
let us take the middle points $(q_0,\dots,q_{N-1})$ of its edges and
transform $\gamma$ into another $N$-gon $T(\gamma)$:
$$
\gamma=(x_0,\dots,x_{N-1}) \stackrel{T}{\longrightarrow}
(q_0,\dots,q_{N-1}) \in \Pi_N.
$$
The deformation $T: \Pi_N \to \Pi_N$ is continuous and decreases the
lengths of curves except geodesics. However the values of the
one-valued functional (\ref{magnetic}) under an analogous
deformation may grow: $S(T(\gamma))
> S(\gamma)$, and in the limit a nontrivial relative cycle $\xi \in
H_\ast(\{S < c+\varepsilon\},\{S < c - \varepsilon\})$ may suspend
on ``a critical point at infinity'' lying on the critical level
$S=c$ (see discussion in \cite{Taimanov}). \footnote{In the case of
such a lack of compactness it is said that $S$ does not satisfy the
Palais--Smale conditions for applying the Morse theory in the
infinite-dimensional situation.}

This difficulty with the convergence of gradient deformations hinders the
application of variational methods to this problem.
The only general statement proved by variational methods for
higher dimensions was obtained in
\cite{Bahri} under some geometrical restrictions on the Riemannian metric and
the magnetic field.

More results were obtained in the two-dimensional situation
\cite{Taimanov,Grinevich}. In particular, the existence of a simple
locally minimal closed magnetic geodesic was proved for ``strong''
magnetic fields on two-dimensional compact manifolds
\cite{Taimanov84,Taimanov91}.

In \cite{CMP} for ``non-strong'' exact magnetic fields
\footnote{A magnetic field $F$ is called exact if the two-form $F$ is exact:
$[F] = 0 \in H^2(M;\R)$.}
on two-dimensional compact surfaces it was proved that a closed nontrivial
extremal always exists. This was done
by using the theory of the Man\'e critical
level and the classical Morse theory (in the situation of Finsler metrics)
and together with the existence result for ``strong'' magnetic fields
this implies the existence of a closed magnetic geodesic on every energy level
for an exact magnetic field on a closed two-dimensional manifold.
\footnote{Again by a combination of different methods
(symplectic and variational) it was recently proved that on almost
all energy levels there exists a closed magnetic geodesic for
an exact magnetic field \cite{Contreras}. One of the main ingredients
of the proof is the
result which states that if an energy level of a magnetic geodesic flow is of
contact type then the (free time) action functional meets
the Palais--Smale condition.}

There is another approach initiated by Arnold \cite{Arnold}
and Kozlov \cite{Kozlov}
to the periodic problem for magnetic geodesics and it is
based on the fixed point theorems from symplectic topology
(we refer to a recent article \cite{Ginzburg} which,
in particular, contains a good survey of the results obtained
via this approach).

A completely new approach to establish the existence of periodic magnetic
geodesics on the two-sphere was very recently proposed in
\cite{Schneider}.

We would like to stress that the theory of periodic magnetic geodesics
is still rich by open problems waiting for their solutions:

a) for the case of exact magnetic field the existence of a closed magnetic
geodesic in general case is still an open problem;

b) for non-exact fields the existence result in general does not hold:
there are no closed magnetic geodesics
on a hyperbolic closed two-surface with its area form as a magnetic field
\cite{Ginzburg96}. However it is
interesting to understand how generic is this phenomenon and under which
conditions periodic trajectories exist.

\end{document}